\documentclass[12pt]{amsart}
\usepackage{amsmath}
\usepackage{graphicx}
\usepackage{amsfonts}
\usepackage{amstext}
\usepackage{amsbsy}
\usepackage{amsopn}
\usepackage{upref}
\usepackage{amsthm}
\usepackage{amsmath}
\usepackage{amssymb}

\newtheorem{thm}{Theorem}[section]
\newtheorem{cor}[thm]{Corollary}
\newtheorem{lemma}[thm]{Lemma}
\newtheorem{prop}[thm]{Proposition}
\theoremstyle{definition}

\theoremstyle{remark}

\numberwithin{equation}{section}

\newcommand{\noi}{\noindent}
\newcommand{\sm}{\smallskip}
\newcommand{\me}{\medskip}
\newcommand{\bi}{\bigskip}
\newcommand{\RR}{{\rm Re}}
\newcommand{\II}{{\rm Im}}

\newcommand{\si}{\sigma}
\newcommand{\de}{\delta}
\newcommand{\Om}{\Omega}
\newcommand{\dOm}{\partial\Omega}
\newcommand{\D}{\mathbb{D}}

\def\dD{\partial{\mathbb D}}

\def\z{{\zeta}}
\def\dis{\displaystyle}

\newcommand{\Ci}{\hat{\mathbb{C}}}
\newcommand{\Cc}{\mathbb{C}}
\newcommand{\N}{\mathcal{N}}
\newcommand{\No}{\mathcal{N}_0}
\newcommand{\Ns}{\mathcal{N}^*}
\newcommand{\C}{\mathcal{C}}
\newcommand{\Co}{\mathcal{C}_0}

\def\af{1+a_2f}
\def\bet{\begin{thm}}
\def\ent{\end{thm}}
\def\bec{\begin{cor}}
\def\enc{\end{cor}}
\def\bel{\begin{lemma}}
\def\enl{\end{lemma}}
\def\bep{\begin{prop}}
\def\enp{\end{prop}}
\def\bp{\begin{proof}}
\def\ep{\end{proof}}
\def\bq{\begin{equation}}
\def\eq{\end{equation}}
\def\ri{\rightarrow}
\def\Ra{\Rightarrow}
\def\fd{f^{*}}
\def\mR{\mathcal{R}}
\def\mP{\mathcal{P}}
\def\mL{\mathcal{L}}
\def\mH{\mathcal{H}}

\begin{document}

\title[Ahlfors-Weill reflection]{The Ahlfors-Weill reflection on convex domains and Nehari quasidisks}
\author{Martin Chuaqui}
\thanks{{\sl Key words: convex domains; Nehari class; Ahlfors-Weill reflection; quasidisk; second coeffi-}
\endgraf {{\sl cient} }
\endgraf {\sl 2000 AMS Subject Classification}: 30C45, 30C35.}
%


\begin{abstract}
The estimate $$\RR\{a_2f\}>-\frac12$$ derived for convex mappings in \cite{FMR}, is interpreted here in terms of the Ahlfors-Weill reflection to show that for such domains $\Om$, the mediatrix of the segment $[w, \mR_w]$ joining a point $w\in\Om$ and its reflection $\mR_w$ lies always outside $\Om$. In particular, the midpoint of the segment is also outside
$\Om$. We determine the extremal cases when such a midpoint can lie of the boundary $\partial\Om$. The normalization $\fd=\frac{f}{1+a_2f}$ to a Möbius equivalent mapping with vanishing second coefficient leads to important distinctions between bounded an unbounded domains. We finally derive a geometric characterization of Nehari quasidisks in terms of the distance to the boundary of the Ahlfors-Weill reflection.

\end{abstract}
\maketitle

\section{Introduction}

The purpose in this paper is to establish new geometric properties of the Ahlfors-Weill extension for mappings in the Nehari class and the sublcass of convex mappings. Much of our motivation stems from the interesting paper \cite{FMR}, where the authors establish the crucial estimate
\bq
\RR\{a_2f\}>-\frac12
\eq
whenever $f$ is a convex mapping of the unit disk $\D$. Here $a_2$ represents as usual the second Taylor coefficient. We interpret this in terms of the  Ahlfors-Weill construction to show that, for convex mappings, the mediatrix of the segment 
$[w,\mR_w]$ joining a point $w\in f(\D)$ and its Ahlfors-Weill reflectrion $\mR_w$, lies always outside the image domain. This is Theorem 3.1. Extremal cases are analyzed when such a line can make contact with the boundary of the domain. 

\sm

The Nehari class will be denoted by $\N$, and consists of the locally univalent functions $f$ of the unit disk $\D$ with $f(0)=0, f'(0)=1$ that satisfy
\bq 
(1-|z|^2)^2|Sf(z)|\leq 2  \, ,
\eq
where $Sf=(f''/f')'-(1/2)(f''/f')^2$ is the Schwarzian derivative. Mappings in $\N$ are univalent in $\D$ \cite{N1} and spherically continuous in $\overline{\D}$ \cite{GP}. They remain univalent in the closed disk except when $f(\D)$ is, up to a Möbius transformation, a parallel strip \cite{GP}. It was shown in \cite{N2} that $\C\subset\N$ for the class $\C$ of convex mappings defined on $\D$. We will often consider the normalization $f''(0)=0$ closely related to whether image domains are bounded or not. The corresponding normalized classes will be denoted by $\No$ and $\Co$. Every function in $\N$ has an analytic counterpart in $\No$ via the Möbius shift
\bq 
\fd=\frac{f}{\af} \, ,
\eq
which preserves both value and derivative at the origin. The mapping $\fd$ is analytic in $\D$ because
\bq 
-\frac{1}{a_2}\notin f(\D) 
\eq
for $f\in\N$ (see \cite{CO1}) and belongs to $\No$ in light of the Möbius invariance of the Schwarzian. Furthermore, $-1/a_2$ lies on the boundary of $f(\D)$ if and only if $f(\D)$ is Möbius conjugate to a parallel strip \cite{CO2}. In \cite{FMR}, the condition (1.4) is rephrased as $-1\notin a_2f(\D)$, and it is shown that for $\Om_f=a_2f(\D)$
\bq
\underset{f\in\N}\cup\Om_f=\Cc\backslash\{-1\}\, .
\eq
As a complement, we show in Lemma 3.2 that for a sequence $a_{2,n}f_n(z_n)$ of values to approach $-1$, a subsequence $f_{n_k}$ must converge to a parallel strip mapping. We also note that the above normalization will not always produce a mapping in $\Co$ when $f\in\C$, a question that will not be analyzed here. See \cite{BS1, BS2} for related issues regarding Möbius transformations of convex domains.

\me

Section 2 will be devoted to convex mappings, for which the much stronger estimate (1.1) holds \cite{FMR}. We will prove that
\bq
\RR\{a_2f\}=-\frac12
\eq 
can only occur on finite boundary points when $f$ maps onto a halfplane or an infinite (convex) sector. These extremal configurations are mentioned in \cite{FMR}, without ruling out other cases. In particular, equality can never hold on the boundary of a bounded convex domain. It is interesting to note that, in addition to unbounded convex domains other than a sector or a halfplane, unbounded non-extremal cases with $a_2\neq 0$ include paralel strip mappings with $f(0)$ not on the axis of symmetry.
In Section 2 we will take a different approach to the procedure of normalizing, showing that for convex mappings, $a_2\fd$ is always bounded by 1, and strictly bounded below 1 precisely when $f$ is itself bounded. Finally, we will establish a new characterization of Nehari quasidisks that does not involve the Schwarzian or the pre-Schwarzian, but instead, that $d(w,\dOm)$ and $d(\mR_w,\Om)$ are comparable near $\dOm$. 

\section{Convex domains}

We recall briefly the proof of (1.1) given in \cite{FMR}. For $f\in \C$ and $\z\in\D$ the function
$$
F(z)=\frac{\z z}{f(\zeta)}\frac{f(z)-f(\z)}{z-\z}
$$
is starlike of order $1/2$, hence 
\bq 
G(z)=\dis z\frac{F'(z)}{F(z)}
\eq 
satisfies
$$
\RR\{G(z)\}>\frac12 \, ,
$$
see  \cite{S-S, Su}. The function $-1+1/G$ is therefore a selfmap of $\D$ and, since it fixes the origin, it has the form 
\bq 
-1+\frac{1}{G(z)}=zh(z)
\eq 
with $|h|\leq 1$.  The Schwarz-Pick lemma applied to $h$ gives
\bq 
|a_2(G)-a_1^2(G)|\leq 1-|a_1(G)|^2 \, ,
\eq
from which the authors derive that
\bq
\RR\{a_2f(\z)\}\geq -\frac12+\frac12(1-|\z|^2)\left|\frac{f(\z)}{\z}\right|^2\geq-\frac12 \, . 
\eq
The same argument can be applied when $|\z|=1$ as long as $f(\z)$ is a finite boundary point. They point out that (2.4) is sharp in the sense that equality will hold for points on $\dD$ when $f$ maps to a halfplane or a sector. We show here that these are the only possible cases for equality in (2.4).

\bet
Let $f\in\C$ and suppose that there exists a finite boundary point $f(\z), \, \z\in\dD,$ for which
\bq 
\RR\{a_2f(\z)\}=-\frac12 \, .
\eq
Then $f$ maps $\D$ either to a halfplane or an infinite sector. 
\ent

\bp 
Since convex mappings are (spherically) continuous in the closed disk, the inequality
$$ \RR\{a_2f\}\geq-\frac12$$ 
will hold on $\dD$. Let $\z\in\dD$ be such that (2.5) holds. After a rotation of $f$, we may assume that $\z=-1$. Both (2.4), and therefore (2.3), must be equalities, hence the function $h$ in (2.2) is either a unimodular constant or an automorphism of the disk. The explicit relation between $f$ and $h$ is derived from (2.1) and (2.2) in the form

\bq
(1+z)(1+zh)f'=(1-h)(f-b) 
\eq
for $b=f(-1)$ .

Suppose first that $h=c$ is a unimodular constant. From (2.6) we see that $c\neq 1$, and a simple integration gives that
$$
f(z)=\frac{z}{1+cz} \, ,
$$
that is, a halplane mapping, for which (2.5) holds for all $\z\neq -\bar{c}$.

\sm

Assume now that $h$ is an automorphism of $\D$, written in the form
$$
h(z)=c\frac{z+a}{1+\bar{a}z}
$$
for some $|a|<1$ and $|c|=1$. By evaluating at $z=0$ both (2.6) and the equation obtained from differentiation, we see that 
\bq
b=\frac{1}{h(0)-1} \, ,
\eq
and
\bq
a_2=-h(0)+\frac12bh'(0)\, .
\eq
Using that $h(0)=ac$ and $h'(0)=c(1-|a|^2)$ we see that 
$$
a_2b=\frac{2ac(1-ac)+c(1-|a|^2)}{2(1-ac)^2} 
$$

$$
\hspace{3.4cm}=-\frac12+\frac{1-a^2c^2+c(1-|a|^2)}{2(1-ac)^2}=-\frac12+\de \, ,
$$

\sm
\noi
hence we must have $\RR\{\delta\}=0$, that is, 
\sm
$$
\RR\{\left(1-a^2c^2+c(1-|a|^2)\right)(1-\bar{a}\bar{c})^2\}=0 \, .
$$

\sm
\noi
The first term in the multiplication of both factors in the real part is given by
\sm
$$
(1-a^2c^2)(1-\bar{a}\bar{c})^2=|1-ac|^2(1+ac)(1-\bar{a}\bar{c})
$$
$$
\hspace{4.9cm}=|1-ac|^2\left(1-|a|^2+2i\II\{ac\}\right) \, ,
$$
and thus 
$$
\RR\{\delta\}=(1-|a|^2)\RR\{|1-ac|^2+c(1-\bar{a}\bar{c})^2\}=0\, .
$$

\sm
\noi
Since $|c(1-\bar{a}\bar{c})^2|=|1-ac|^2$ we conclude that $c(1-\bar{a}\bar{c})^2=-|1-ac|^2$, which gives
\bq
c=-\frac{1-\bar{a}}{1-a}\, .
\eq
We observe that $c$ can never be equal to $1$.
Using this in (2.6) we obtain
\bq
\frac{f'}{f-b}=\frac{1-|a|^2}{(1-a)(cz+1)(z+1)} 
 = \beta\left(\frac{1}{z+1}-\frac{c}{cz+1}\right) \, ,
\eq
where
$$
\beta=\frac{1-|a|^2}{2(1-\RR\{a\})}\in (0,1) \, .
$$

\sm
\noi
which upon integration, gives
$$
f(z)=-b\left[\left(\frac{1+z}{1+cz}\right)^{\beta}-1\right] \, .
$$
Since $c\neq 1$, this function corresponds to a mapping onto a convex sector with opening angle equal to $\beta$. Equality in (2.5) will indeed happen for $\z=-1$, while strict inequality will hold for all  
$\z$ not equal to $-1$. For real values of $a\in(-1,1)$, this function becomes the more familiar mapping onto such a sector given by

$$
\frac{1}{2a}\left[\left(\frac{1+z}{1-z}\right)^a-1\right] \, .
$$

\ep

\bec
If $f\in\C$ is bounded then 
\bq 
\RR\{a_2f(z)\}\geq b>-\dis\frac12
\eq
for all $z\in\D$ and some constant $b$.
\enc

\bp
The function $\RR\{a_2f(\z)\}$ is continuous for $\z\in\dD$ and, by Theorem 2.1, it cannot reach the value $-1/2$ because $f$ is bounded. This proves the theorem.

\ep

\me
Recall that the only unbounded convex mapping (or even Nehari mapping) with vanishing second coefficient is 
$$ L(z)=\dis\frac12\log\frac{1+z}{1-z} \, ,
$$ 
with $L(\D)$ equal to a parallel strip \cite{CO1}. The other Möbius equivalent presentations of a parallel strip are two circles tangent at a finte point or a circle tangent to a line at a finte point, but they fail to be convex. Nevertheless, mappings onto a parallel strip with $a_2\neq0$ can be obtained by suitable Koebe transforms of $L$ that map $0\in\D$ to a point not on the axis of symmetry of the strip. We will refer to them as non-symmetric parallel strip mappings. 

\bet
Let $f$ be a conformal map onto an unbounded convex domain with $a_2\neq 0$. Then
\bq
|a_2\fd(z)|<1 \; , \; \; z\in\D
\eq
and $|a_2\fd(\z)|=1$ at some $\z\in\dD$. More precisely, 

\bi
\noi
$(i)$ if $f$ is not a halfplane or a sector mapping, then $|a_2\fd(\z)|=1$ will hold at a unique point on $\dD$ corresponding to the point at $\infty$ on $\partial f(\D)$;

\me
\noi
$(ii)$ for a non-symmetric parallel strip maping, there will be two solutions corresponding to each end at $\infty$; 

\me
\noi
$(iii)$ for a sector mapping, there will be two solutions corresponding to the vertex and the point at $\infty$; 

\me
\noi
$(iv)$ and for a halfplane, the equation will hold on the entire unit circle.  
\ent

\bp
The relation
\bq
a_2\fd=\frac{a_2f}{1+a_2f}
\eq 

\me
\noi
between $f$ and its normalization $\fd$ can be expressed in terms of the mapping
$$
w=\frac{z}{1+z} \, ,
$$
with $\RR\{z\}>\dis-\frac12$ corresponding to $|w|<1$. This proves (2.10). 

For the remaining assertions, simply observe that $|a_2\fd(\z)|=1$ either when $\RR\{a_2f(\z)\}=-1/2$ at a finite boundary point or when $f(\z)=\infty$. 

\ep

\bet
Let $f$ be a conformal map onto a bounded convex domain. Then
$$
|a_2\fd|\leq c<1 
$$
for some constant $c$.
\ent

\bp
Since $f$ is bounded, the function $\RR\{a_2f(\z)\}$ is continuous on $\dD$ and must be strictly larger then $-1/2$ by Theorem 2.1. This proves the theorem. 

\ep

\section{The Ahlfors-Weill reflection and Nehari quasidisks}

\bi

The well-known Ahlfors-Weill extension for mappings $f$ with 
$$
(1-|z|^2)^2|Sf(z)|\leq 2t<2
$$ 
can be viewed as a reflection across the boundary of the image $\dOm$ having desired quasiconformal properties. This reflection is given by

$$
f(z)\ri f(z)+\frac{(1-|z|^2)\dis f'(z)}{\dis\bar{z}-(1-|z|^2)(f''/2f')(z)} \, ,
$$

\me
\noi
and provides a homeomorphic extension for all $f\in\N$ not Möbius conjugate to $L$ \cite{CO2}. In this larger class, the extension will not be necessarily quaisiconformal, for example, when $f(\D)$ is an unbounded convex domain with a cusp at infinity \cite{CDO}. For a convex domain $\Om$, let $\mR:\Om\ri\Ci$ be the Ahlfors-Weill reflection, and denote $\mL_w$ the mediatrix of the line segment $[w,\mR_w]$ with midpoint $\mP_w$. We will denote by $\mH_w$ the halfplane determined by $\mL_w$ that contains the point $w$.

\bet
For a convex domain $\Om$ and all $w\in\Om$ we have
\bq 
\Om\hspace{.02cm} \cap \mL_w=\emptyset \, .
\eq 
In particular, $\mP_w\notin\Om$.
\ent

\bp
The proof is based on linear invariance and the following simple geometric description of the set of points $w$ for which $\RR\{aw\}>-1/2$. If $a>0$ then the set corresponds to the halfplane 
$\RR\{w\}>-1/(2a)$. In general, and after a rotation, we see the set corresponds to the halfplane containing the origin determined by the line orthogonal to the segment $[0,-1/a]$ that passes through its midpoint. Because $-1/a_2$ is the Ahlfors-Weill reflection of $w=0$, we see that (3.1) holds for $w=0$ by setting $a=a_2$.

For a generic point $z_0\in\D$ we consider the Koebe transform
\bq
g(z)=\frac{f(\si(z))-f(z_0)}{(1-|z_0|^2)f'(z_0)}
\eq
with $\si(z)=\dis\frac{z+z_0}{1+\bar{z_0}z}$. Then $g\in\C$ with a second coefficient given by
$$
b_2=(1-|z_0|^2)\frac{f''(z_0)}{2f'(z_0)}-\bar{z_0} \, ,
$$
and therefore $\RR\{b_2g(z)\}>-1/2$. We conclude that $g(\D)$ is a subset of the halfplane containg the origin determined by the mediatrix of the segment $[0,-1/b_2]$. 
The conclusion that $f(\D)\subset\mH_{f(z_0)}$ follows after applying the affine conformal transformation $w\ri Aw+B$ with $A=(1-|z_0|^2)f'(z_0)$ and 
$B=f(z_0)$.
\ep

Recall that for a mapping $f\in\Ns$, the value $-1/a_2$ lies outside the closure of the image $f(\D)$. For $a_2\neq 0$, let
\bq
\de_f=\underset{z\in\D}\inf \, \left|f(z)+\frac{1}{a_2}\right| \, ,
\eq
which is to be understood in terms of the spherical distance  when $a_2=0$.

\bel
Let $\{f_n\}\subset\Ns$ be a sequence of Nehari mappings with 
\bq
\lim_{n\ri \infty}\de_{f_n}\, = \, 0 \, .
\eq
Then some subsequence of  normalized mappings $f_{n_k}^*$ converges locally uniformly to $L$.
\enl

\bp
Since $f_n\in\Ns$, each normalized mapping $f_n^*$ is bounded by some constant $C_n<\infty$. Because of the Möbius invariance of the Ahlfors-Weill reflection, the condition (3.4) implies that the spherical distance between the domains $f_n^*(\D)$ and $\infty$ tends to 0. Hence $C_n\ri \infty$, and the conclusion follows from \cite [Corollary 4, (ii) ]{JQN}.
\ep

\bet
Let $\Om=f(\D)$ for a mapping $f\in\No$. Then the following are equivalent:

\me
\noi 
(i) $\Om$ is a quasidisk;

\me
\noi
(ii) there exists $c>0$ such that $a_2(g)g$ omits a fixed neighborhood $|w+1|<c$ of 

\noi
the point $-1$ for every Koebe transform $g$ of $f$;

\me
\noi 
(iii) there exists $c>0$ such that $d(\mR_w,\Om)\geq cd(w,\dOm)$ for all points $w\in\Om$.
\ent

\bp
To show $(i) \Ra (ii)$ via the contrapositive, we assume that $\de_{g_n}\ri 0$ for some sequence of Koebe transforms of $f$. By Lemma 3.2, some subsequence of normalized mappings $g_{n_k}^*$ converges to $L$ which, in turn, implies that some subsequence of $g_{n_k}$ (still denoted by $g_{n_k}$) converges to a mappings of the form $L_1=L/(1+aL)$. It can be verified that the Koebe transforms of $L_1$ as in (3.2) with
$$\si(z)=\frac{z+x}{1+xz} \;\, , \; \;x\in(0,1)$$
converge to $L$ as $x\ri 1$. Hence, approriate Koebe transforms of of the $g_{n_k}$ will converge to $L$, which implies that $\Om$ cannot be a quasidisk by \cite[Theorem 5]{JQN}.

\me

\noi
$(ii)\Ra(iii)$. We repeat the proof of Theorem 3.1, except that missing the halplane $\RR\{a_2(g)g\}<-\frac12$ is replaced is by missing the disk $|w+1|<c$. The proof shows that the constant $c$ is the same in both parts.

\me

\noi
$(iii)\Ra(i)$. Again, by contrapositive, if $\Om$ were not a quasidisk, then by \cite[Theorem 5]{JQN} a sequence of Koebe transforms $g_n$ would converge to $L$. Then $a_2(g_n)\ri a_2(L)=0$, and therefore, the normalized sequence $g_n^*$ would also converge to $L$ locally unformly. By Lemma 3.2, $(iii)$ could not hold. This finishes the proof.
\ep

\bibliographystyle{plain}

\bi
\noi
{\small Facultad de Matemáticas, Pontificia Universidad Católica de Chile,
Santiago, Chile,\, \email{mchuaqui@uc.cl}

\end{document}